\documentclass{article}
\usepackage{amsmath}
\usepackage{amssymb}
\usepackage{latexsym}
\usepackage{amsthm}
\usepackage{mathrsfs} 
\usepackage{wasysym}
\usepackage{fancyhdr}
\usepackage{amsfonts}
\usepackage{amssymb}
\usepackage{epsfig}
\usepackage{epstopdf}
\usepackage[normalem]{ulem}
\usepackage{eucal}
\usepackage{tikz}

\newtheorem{thm}{Theorem}[section]
\newtheorem{cor}[thm]{Corollary}

\newtheorem{prop}[thm]{Proposition}

\newtheorem{obs}[thm]{Observation}

\begin{document}

\title{On Sign-Invertible Graphs}

\author{Isaiah Osborne\thanks{Department of Mathematical Sciences, Middle Tennessee State University, Murfreesboro, TN 37132. Email: imo2d@mtmail.mtsu.edu. Partially supported by a URECA grant of MTSU.} \, and Dong Ye\thanks{Department of Mathematical Sciences and Center for Computational \& Data Science, Middle Tennessee State University, Murfreesboro, TN 37132. Email: dong.ye@mtsu.edu. Partially supported by a grant from Simons Foundation (359516).}}
\date{December 22 2022}

\maketitle

\begin{abstract}

Let $G$ be a graph and $A$ be its adjacency matrix. A graph $G$ is invertible if its adjacency matrix $A$ is invertible and the inverse of $G$ is a weighted graph with adjacency matrix $A^{-1}$. A signed graph $(G,\sigma)$ is a weighted graph with a special weight function $\sigma: E(G)\to \{-1,1\}$. A graph is sign-invertible (or sign-invertible) if its inverse is a signed graph. 
A sign-invertible graph is always unimodular. The inverses of graphs have interesting combinatorial interests. 
In this paper, we study inverses of graphs and provide a combinatorial description for sign-invertible graphs, which provides a tool to characterize sign-invertible graphs. As applications, we complete characterize sign-invertible bipartite graphs with a unique perfect matching, and sign-invertible graphs with cycle rank at most two. As corollaries of these characterizations, some early results on trees (Buckley, Doty and Harary in 1982)  and unicyclic graphs with a unique perfect matching 
(
Kalita and Sarma in 2022) follow directly.


\medskip

 \medskip


\end{abstract}

\section{Introduction}

All graphs considered in this paper do not have multi-edges but may have loops. Let $G$ be a graph with vertex set $V(G)$ and edge set $E(G)$.  If $G$ has $n$ vertices, its 
adjacent matrix of $G$ is an $n\times n$ symmetric $\{0,1\}$-matrix, denoted by $A (G)=[a_{ij}]$ such that, for each pair of vertices $i$ and $j$, the $(i,j)$-entry
$a_{ij}=1$ if $ij\in E(G)$ and $a_{ij}=0$ otherwise. Particularly, $a_{ii}=1$ if $G$ has a loop at the vertex $i$. The eigenvalues of $A(G)$ is also called the eigenvalues of the graph $G$, which are usually ordered in a non-increasing order such as $\lambda_1(G) \ge \lambda_2(G)\ge \cdots \ge\lambda_n(G)$ where $n=|V(G)|$.


A graph $G$ is  {\em invertible} or {\em non-singular} if its adjacency matrix
is invertible \cite{YYMK}, equivalently has a non-zero determinant. If $G$ is invertible, then the inverse of its adjacency matrix is a symmetric matrix which is an adjacency matrix of a weighted graph which may have loops.
The inverse of a graph $G$ could be represented as a weighted graph, denoted by $(G^{-1},w)$, such that $V(G^{-1})=V(G)$,  two vertices $i$ and $j$ are adjacent if  the $(i,j)$-entry of $A(G)^{-1}$ is not zero, and the weight function $w(ij) =(A(G)^{-1})_{ij} \ne 0$. If $w(ii)\ne 0$, then the inverse of $G$ has a loop at vertex $i$. A graph $G$ is {\em sign-invertible} or {\em sign-invertible} if its inverse $(G^{-1},w)$ has a weight function $w: E(G^{-1})\to \{-1,+1\}$  (ref.~\cite{BDH}), and such weighted graph is also called a {\em signed graph}.  If a graph $G$ is invertible, then its smallest positive eigenvalue is equal to the reciprocal of the largest eigenvalue of its inverse, and its largest negative eigenvalue is equal to the reciprocal of the smallest eigenvalue of its inverse. 

Graph inverse has been studied since 1970s (cf. \cite{CD85, GR79, K}). A particular attention of graph inverse has been paid to bipartite graphs with a unique perfect matching~\cite{CD85, TK}. It is known that the inverse of the adjacency matrix of a bipartite graph with a unique perfect matching is an integral matrix~\cite{CD85, TK, YYMK}. There is a correspondence relation between bipartite graphs with a unique perfect matching and posets (see~\cite{CD85, YY18}). The inverses of bipartite graphs with a unique perfect matchings have many combinatorial interests such as M\"{o}bius function of posets~\cite{CD85,YY18}, Motzkin numbers~\cite{MM} and etc. Another motivation for study graph inverse comes from the problem to bond median eigenvalues~\cite{CGS,YYMK}. Median eigenvalues of a graph have physical meanings~\cite{GR79} and the absolute value of their differences is called the HOMO-LUMO gap, an important parameter which could be used to estimate molecule's chemical stability~\cite{GR79, Mo, WYY18}.  If a graph is bipartite and invertible, the two median eigenvalues are the smallest positive eigenvalue and the largest negative eigenvalue of the graph, which can be estimated by investigating the largest and the smallest eigenvalues of its inverse. Using graph inverse, the median eigenvalues of bipartite graphs with a unique perfect matchings~\cite{WYY18}, and stellated graphs of trees, and corona graphs~\cite{YYMK} have been proven belonging to the interval $[-1,1]$. For a given family of graphs, it is an interesting problem to ask which graphs maximize the HOMO-LUMO gap. The answer to this problem is only known for trees~\cite{SH}, and is open even for unicyclic graphs, and bipartite graphs with a unique perfect matching.

Graph inverse is extensively studied for trees, unicyclic graphs and bipartite graphs with a unique perfect matching.  Harary and Minc~\cite{FH76} show that there is exactly one connected graph with a inverse being a graph which is $K_2$, a complete graph on two vertices. All trees with a perfect matching are sign-invertible~\cite{BDH}, i.e. having a signed graph as inverse. A signed graph $(G,\sigma)$ is balanced if every cycle of $(G,\sigma)$ has an even number of negative edges (ref~\cite{TZ}).  It is known that a signed graph $(G,\sigma)$ is balanced if and only if it has an edge-cut such that an edge is negative if and only if it belongs to the edge-cut~\cite{FH53}. In~\cite{CD85}, Godsil asked  which bipartite graphs with
a unique perfect matching a balanced signed inverse, and settled the problem for all trees.  A partial result of the problem was obtained by Akbari and Kirkland~\cite{AK} for unicyclic graphs. This problem was completely settled by Yang and the second author in~\cite{YY18} based on the combinatorial description of graph inverse~\cite{BNP, YYMK}. The problem on weighted graphs are also studied for positive weight functions and some partial results on trees and unicyclic graphs have been obtained~\cite{PP15, PP17}. The inverses of trees and unicyclic graphs have been studied in terms of the reciprocal eigenvalue properties~\cite{BPP, BNP, PP15}. Bipartite unicyclic graphs whose inverses have underlying graphs being unicyclic or bicyclic are also studied~\cite{SP17, TK}. Akabari and Kirkland~\cite{AK} characterized sign-invertible bipartite unicyclic graphs with a unique perfect matching. Recently, a characterization of sign-invertible non-bipartite unicyclic graphs with a unique perfect matching has been obtained by Kalita and Sarma~\cite{KS}. Even though many efforts have been made to study sign-invertibility of graphs, a complete characterization remains unknown for unicyclic graphs.

In this paper, we study sign-invertibility of graphs and the goal is to develop tools to study the problem. A graph $G$ is {\em unimodular} if its adjacency matrix has determinant either $-1$ or $+1$. Unimodularity is a 
necessary condition for a graph to be sign-invertible.  We proved that a graph $G$ is sign-invertible if and only if it is unimodular and each pair of vertices has to satisfies some combinatorial properties. The details of the combinatorial description will be presented in Section 3. Based on this result, we are able to characterize the sign-invertibility for bipartite graphs with a unique perfect matching and  graphs with cycle rank at most two. As direct corollaries of these results, previous results on trees~\cite{BDH} and unicyclic graphs with a unique perfect matching~\cite{AK, KS} follow.


\section{Determinant and unimodularity }

A {\it 2-matching} of a graph $G$ is a spanning subgraph of which every vertex has degree at most two. A 2-matching is {\em perfect} if every component of the 2-matching is either a cycle or a single edge (or a complete graph with two vertices $K_2$).  A perfect 2-matching is also called an elementary subgraph or a Sach's subgraph in literature. A {\em perfect matching} is a    perfect 2-matching without cycles. A well-known result of Edmonds~\cite{Ed} show that it is polynomial time to find a perfect matching in a given graph with perfect matchings, which also implies that it is polynomial time to find a perfect 2-matching in a graphs with perfect 2-matchings due to Tutte's reduction~\cite{Tutte}.  However, Papadimitriou~\cite{CP} proved that it is NP-complete to determine whether a graph has a perfect 2-matching without cycles of length five or less.
For a perfect 2-matching $H$ of a graph $G$, we always denote it  by  $H = C_H\cup M_H$,
where $C_H$ consists of all cycles of $H$, and $M_H$ consists of disjoint edges. The edge set $M_H$ is a matching of $G$, where a matching is a set of disjoint edges.

\vskip 2mm
For convenience, the determinant of adjacency matrix of a graph $G$ is denoted by $\det(G)$. If $G$ is an empty graph (having no vertices), for convenience, set $\det(G)=1$. The following classic result shows the determinant of the adjacent matrix of a graph can be computed via perfect 2-matchings.

\begin{thm}[Harary, \cite{FH62}]\label{determinant}
Let $G$ be a simple graph. Then
$$
\det(G)=\sum\limits_{H}2^{|C_H|}(-1)^{|C_H|+|E(H)|}
$$
where $H = C_H\cup M_H$ is a perfect 2-matching of $G$.
\end{thm}

By the above theorem, the edges of a graph $G$  not belong to any perfect 2-matching do not affect the determinant of the graph. Let $E_0$ be the set of edges of $G$ which does not belong to a perfect 2-matching. Then $\det(G)=\det(G-E_0)$. It follows from a property of block matrices, we have the following proposition. 

\begin{prop}\label{prop:det}
Let $G$ be a graph and let $E_0$ be the set of edges of $G$ which are not contained by any perfect 2-matchings. Assume $G_1, G_2, \ldots, G_k$ are components of $G-E_0$. Then
\[\det(G)=\prod_{i=1}^k \det(G_i).\]
\end{prop}

By Proposition~\ref{prop:det} and Theorem~\ref{determinant}, the following results hold trivially. 

\begin{prop}\label{prop:forest}
Let $G$ be a graph without cycles. Then $|\det(G)|=1$ if and only if $G$ has a perfect matching.
\end{prop}

A direct application of Theorem~\ref{determinant} on cycles will give the following proposition.

\begin{prop}[Akbari and Kirkland \cite{AK}]\label{prop:cycle}
Let $C_n$ be a cycle of $n$ vertices. Then 
\begin{equation*}
\det(C_n)=
\left \{ \begin{array}{lll}
0 & \mbox{ if } n\equiv 0 \pmod 4;\\
-4 & \mbox{ if } n\equiv 2 \pmod 4;\\
2 & \mbox{ if } n\equiv 1\pmod 2.
\end{array}
\right.
\end{equation*}
\end{prop}
 
A graph $G$ is {\em unimodular} if $|\det(G)|=1$. The above proposition shows that a cycle is not unimodular. The determinants of cycles seems easy to calculate. But for graphs containing more than one cycles, it becomes complicated. For examples, the following two propositions deal with two families of graphs which are not far from a cycle. 

A {\em theta graph} consists of three internally disjoint paths with $\theta_i$ vertices for $i\in [3]$ joining two vertices $x$ and $y$, which is denoted by $\Theta(\theta_1,\theta_2,\theta_3)$. The two vertices $x$ and $y$ are called {\em central vertices} and the three paths joining $x$ and $y$ are called {\em central paths}. Note that, the ordering of $\theta_1, \theta_2, \theta_3$ in the notation does not matter. For convenience, 
we sometimes assume that $2\le \theta_1\le \theta_2\le \theta_3$. For example, see Figure~\ref{fig:theta}. \medskip

\begin{figure}[htbp]
\begin{center}
\begin{tikzpicture}
\draw[](.5,.5)--(2,2)--(4,2)--(5.5,.5);
\draw[] (.5,.5)--(5.5,.5);
\draw[](.5,.5)--(1.5,-.7)--(3,-1)--(4.5,-.7)--(5.5,.5); 
\fill (2,2) circle (2pt);
\fill (4,2) circle (2pt); 
\fill (0.5,0.5) circle (2pt); 
\fill (5.5, .5) circle (2pt);
\fill (3,.5) circle (2pt);
\fill (4.5,-.7) circle (2pt);  
\fill (1.5,-.7) circle (2pt); 
\fill (3,-1) circle (2pt); 
\end{tikzpicture}
\caption{A theta graph $\Theta(3, 4, 5)$.}\label{fig:theta}
\end{center}
\end{figure}
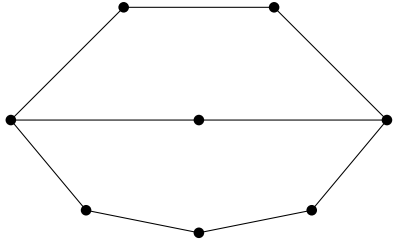

\begin{prop}\label{prop:theta} 
A theta graph $\Theta(\theta_1,\theta_2, \theta_3)$ is unimodular if and only if $\theta_1, \theta_2$ ad $ \theta_3$ are even but not all congruent modulo $4$. 
Furthermore, $\det (\Theta(\theta_1,\theta_2,\theta_3))= 0$ if and only if one of the following holds:\\
{\upshape (i)} $ \theta_1\equiv \theta_2\equiv \theta_3\equiv 1 \pmod 2$, or\\
{\upshape (ii)} exactly one $i\in [3]$ such that $\theta_i$ is even and for the others $j, k\in [3]$, $ \theta_j+ \theta_k\not\equiv 0 \pmod 4$, or \\
{\upshape (iii)} exactly one $i\in [3]$ such that $\theta_i$ is odd and for the others $j,k\in [3]$, $ \theta_j\not\equiv \theta_k \pmod 4$.
\end{prop}
\begin{proof} Let $x$ and $y$ be the central vertices, and $P_i$ be the central path of $\Theta(\theta_1, \theta_2, \theta_3)$ with $\theta_i$ vertices for $i\in [3]$, and let $C_{ij}=P_i\cup P_j$ for $i,j\in [3]$. 
By Theorem~\ref{determinant}, 
\begin{equation}\label{eq1}
\det(\Theta(\theta_1, \theta_2, \theta_3))=\sum_{H} 2^{| C_H|} (-1)^{|C_H|+|E(H)|},
\end{equation} where $H= C_H\cup M_H$ is a perfect 2-matching of $\Theta(\theta_1, \theta_2, \theta_3)$. A perfect 2-matching $H$ of a theta graph is either a perfect matching or has an exactly one cycle. \medskip

\noindent{\bf Claim~1.} {\sl A theta graph $\Theta(\theta_1, \theta_2, \theta_3)$  has a perfect 2-matching if and only if at least one of $\theta_i$ is even.} \medskip

\noindent{\it Proof of Claim~1.} If one of $\theta_1, \theta_2$ and $\theta_3$ is even, without loss of generality, assume $\theta_i$ is even. Then $P_i-\{x,y\}$ is a path with an even number $\theta_i-2$ vertices, which has a perfect matching $M_i$. So $\Theta(\theta_1, \theta_2, \theta_3)$ has a perfect 2-matching $C_{jk}\cup M_i$. 

Now, assume that $\Theta(\theta_1, \theta_2, \theta_3)$ has a perfect 2-matching. Suppose to the contrary that all $\theta_1, \theta_2$ and $\theta_3$ are odd. Then $\Theta(\theta_1, \theta_2, \theta_3)$  has $\theta_1+\theta_2+\theta_3-4$ vertices, which is an odd number. So $\Theta(\theta_1, \theta_2, \theta_3)$ has no perfect matching. Then a perfect 2-matching $H$ contains exactly one cycle, say $C$. Then $\Theta(\theta_1, \theta_2, \theta_3)- V(C)$ is an odd path, which contradicts that  $\Theta(\theta_1, \theta_2, \theta_3)- V(C)$ has a perfect matching $M_H$. This completes the proof of Claim~1. \medskip

By Claim~1, $\det(\Theta(\theta_1, \theta_2, \theta_3))=0$ if all $\theta_1, \theta_2$ and $\theta_3$ are odd. In the following, assume at least one of $\theta_1, \theta_2$ and $\theta_3$ is even.  \bigskip

\noindent{\bf Case~1.} Exactly one of $\theta_1, \theta_2$ and $\theta_3$ is even. Without loss of generality, assume $\theta_i\equiv 0\pmod 2$ and $\theta_j\equiv \theta_k\equiv 1\pmod 2$. Then a perfect 2-matching $H$ of $\Theta(\theta_1, \theta_2, \theta_3)$ contains a perfect matching $M_i$ of $P_i-\{x,y\}$. Furthermore, $H-M_i$ is either a cycle $C_{jk}$ or a perfect matching of $C_{jk}$. Let $M_{jk}$ and $M_{jk}'$ be two disjoint perfect matchings of $C_{jk}$. So $\Theta(\theta_1, \theta_2, \theta_3)$ has exactly three perfect 2-matchings, $C_{jk}\cup M_i$, $M_{jk}\cup M_i$ and $M_{jk}'\cup M_i$. Therefore, it follows from Theorem~\ref{determinant} that
\[\setlength{\jot}{5pt} 
\begin{aligned}
\displaystyle \det(\Theta(\theta_1, \theta_2, \theta_3))&= 2 (-1)^{1+|E(C_{jk})|+ |M_i|}+ 2 (-1)^{|M_{jk}|+|M_i|} \\
&= \displaystyle 2(-1)^{1+(\theta_j+\theta_k-2)+(\theta_i-2)/2}+2(-1)^{(\theta_j+\theta_k-2)/2+(\theta_i-2)/2}\\
&=\displaystyle 2 (-1)^{(\theta_i+\theta_j+\theta_k)/2}((-1)^{(\theta_j+\theta_k)/2}+1). 
\end{aligned}
\]
So \[ \det (\Theta(\theta_1, \theta_2, \theta_3))\equiv 0\pmod 2,\] and $\det (\Theta(\theta_1, \theta_2, \theta_3))=0$ if and only if 
\[(\theta_j+\theta_k)/2  \equiv 1 \pmod 2.\]
Therefore, in this case, $\det(\Theta(\theta_1, \theta_2, \theta_3))=0$ if and only if 
$\theta_j+\theta_k\not\equiv 0 \pmod 4$.  

\medskip
\noindent{\bf Case~2.} Exactly two of $\theta_1, \theta_2$ and $\theta_3$ are even. Without loss of generality, assume that $\theta_i$ is odd, and  $\theta_j, \theta_k$ are even. Let $M_j$ and $M_k$ be the perfect matchings of $P_j-\{x,y\}$ and $P_k-\{x,y\}$ respectively. Note that $\Theta(\theta_1, \theta_2, \theta_3)$ has no perfect matching, and it  has exactly two perfect 2-matchings, $C_{ij}\cup M_k$ and $C_{ik}\cup M_j$, where both $C_{ij}$ and  $C_{ik}$ are odd cycles. By Theorem~\ref{determinant}, 
\begin{equation*}
\setlength{\jot}{5pt} 
\begin{aligned}
\det(\Theta(\theta_1, \theta_2, \theta_3))&= 2 (-1)^{1+|E(C_{ij})|+ |M_k|}+ 2 (-1)^{1+|E(C_{ik})|+ |M_j|}\\
&=2(-1)^{(\theta_k-2)/2}+2(-1)^{(\theta_j-2)/2}.
\end{aligned}
\end{equation*}
Then,
\[\det(\Theta(\theta_1, \theta_2, \theta_3))\equiv 0\pmod 2,\]
and 
$\det(\Theta(\theta_1, \theta_2, \theta_3))=0$ if and only if $\theta_j\not\equiv \theta_k\pmod 4$. \medskip

\noindent{\bf Case~3.} All $\theta_1, \theta_2$ and $\theta_3$ are even. Let $C_{ij}=P_i\cup P_j$ for $i,j\in [3]$, and let $M_i$ be the perfect matching of $P_i-\{x,y\}$ for $i\in [3]$. Then $C_{ij}\cup M_k$ is a perfect 2-matching of  $\Theta(\theta_1, \theta_2, \theta_3)$ for distinct $i,j,k\in [3]$. 
Then $\Theta(\theta_1, \theta_2, \theta_3)$ has exactly six distinct perfect 2-matchings: three of them contain exactly one cycle and the remaining three are perfect matchings. 
By Theorem~\ref{determinant},
 \[\setlength{\jot}{5pt} 
\begin{aligned}
\displaystyle \det(\Theta(\theta_1, \theta_2, \theta_3))&=  \sum_{i,j,k\in [3]} 2 (-1)^{1+|E(C_{ij})|+ |M_k|} + 3 (-1)^{(\theta_1+\theta_2+\theta_3-4)/2} \\
&= \displaystyle \sum_{i,j,k\in [3]} 2(-1)^{1+(\theta_i+\theta_j-2)+(\theta_k-2)/2}+3(-1)^{(\theta_i+\theta_j +\theta_k)/2}\\
&=\displaystyle  2 \sum_{k\in [3]}(-1)^{\theta_k/2}  + 3(-1)^{(\theta_i+\theta_j+\theta_k)/2}.
\end{aligned}
\]
It follows immediately that
 $\det( \Theta(\theta_1, \theta_2, \theta_3))\ne 0$. Furthermore, $|\det(\Theta(\theta_1, \theta_2, \theta_3))|=1$ if and only if $\theta_1, \theta_2$ and $\theta_3$ are not congruent modulo 4. 
This completes the proof.
\end{proof}

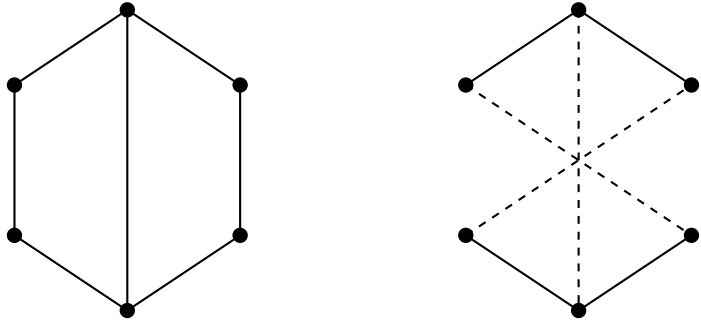
\begin{figure}[htbp]
\begin{center}
\begin{tikzpicture}[thick]

\draw[] (1,2)--(-.5,1)--(-.5,-1)--(1,-2)--(2.5,-1)--(2.5,1)--(1,2);
\draw[] (1,2)--(1,-2);

\filldraw[] (1,2) circle (2.5pt);
\filldraw[] (-.5,1) circle (2.5pt);
\filldraw[] (-.5,-1) circle (2.5pt);
\filldraw[] (1,-2) circle (2.5pt);
\filldraw[] (2.5,1) circle (2.5pt);
\filldraw[] (2.5,-1) circle (2.5pt);

\draw[] (8.5,1)--(7,2)--(5.5,1); \draw[] (8.5,-1)--(7,-2)--(5.5,-1);
\draw[dashed] (7,2)-- (7,-2); \draw [dashed](5.5,1) -- (8.5,-1); \draw[dashed] (5.5, -1)--(8.5,1);
 
\filldraw[] (7,2) circle (2.5pt);
\filldraw[] (5.5,1) circle (2.5pt);
\filldraw[] (5.5,-1) circle (2.5pt);
\filldraw[] (7,-2) circle (2.5pt);
\filldraw[] (8.5,1) circle (2.5pt);
\filldraw[] (8.5,-1) circle (2.5pt);

\end{tikzpicture}
\caption{A theta graph $\Theta(2, 4, 4)$ and its s-inverse (dashed lines are negative edges).}\label{fig:theta}
\end{center}
\end{figure}

A {\em barbell} is a graph consisting of two disjoint cycles which are joined by a path that is internally disjoint from the two cycles. Assume the two cycles have lengths $\theta_1$ and $\theta_2$, and the path has $\theta$ vertices. Then the barbell is denoted by $\mathbf B(\theta_1,\theta_2;\theta)$. A barbell with $\theta=1$ is also called a {\em bowtie} or {\em short barbell}. A vertex of degree bigger than two is called a {\em central vertex} of $\mathbf B(\theta_1,\theta_2;\theta)$ and the path joining the two cycles is called the {\em central path}.

\begin{prop}\label{prop:barbell}
A barbell $\mathbf B(\theta_1,\theta_2;\theta)$ is not unimodular, and 
$\det(\mathbf B(\theta_1,\theta_2;\theta))=0$ if and only if one of the following holds: \\
{\upshape (i)}  $\theta \equiv 1 \pmod 2$ and $\theta_1\equiv \theta_2\equiv 0 \pmod 2$;\\
{\upshape (ii)}  $\theta \equiv 1 \pmod 2$ and $\theta_1 \equiv \theta_2\equiv 1 \pmod 2$, $\theta_1\not\equiv \theta_2 \pmod 4$;\\
{\upshape (iii)} $\theta\equiv 1 \pmod 2$  and  $\theta_1\not \equiv  \theta_2 \pmod 2$ and the even value of $\theta_1$ and $\theta_2$ is congruent to zero modulo 4;\\
{\upshape (iv)} $\theta\equiv 0\pmod 2$ and $\theta_1\equiv \theta_2\equiv 0 \pmod 2$ and they are not both congruent to 2 modulo 4.
\end{prop}
\begin{proof}
Let $C_1$ and $C_2$ be the two cycles of $B(\theta_1,\theta_2;\theta)$ and $P$ be the central path such that $|V(C_i)|=\theta_i$ for $i\in [2]$ and $|V(P)|=\theta$. 

First, we consider the case that $\theta$ is odd. In this case, $\mathbf B(\theta_1,\theta_2;\theta)$ has a perfect 2-matching if and only if at least one of $\theta_1$ and $\theta_2$ is odd. In other words, $\det(\mathbf B(\theta_1,\theta_2;\theta))=0$ if $\theta_1\equiv \theta_2\equiv 0 \pmod 2$. 

If $\theta_1\equiv \theta_2\equiv 1\pmod 2$, then $\mathbf B(\theta_1,\theta_2;\theta)$ has two perfect 2-matchings, both containing exactly one cycle (either $C_1$ or $C_2$). It follows from Theorem~\ref{determinant} that
\[\setlength{\jot}{5pt} 
\begin{aligned}
\det(\mathbf B(\theta_1,\theta_2;\theta))&= 2 (-1)^{1+|E(C_1)|+(|V(C_2)|+|V(P)|-2)/2} + 2 (-1)^{1+|E(C_2)|+(|V(C_1)|+|V(P)|-2)/2}\\
&= -2((-1)^{(\theta_2+\theta)/2} +(-1)^{(\theta_1+\theta)/2}).
\end{aligned}\]
It follows that $\det(\mathbf B(\theta_1,\theta_2;\theta))\equiv 0 \pmod 2$, and furthermore $\det(\mathbf B(\theta_1,\theta_2;\theta))=0$ if and only if $\theta_1\not\equiv \theta_2 \pmod 4$.

Now, consider the case that $\theta_1\not\equiv \theta_2\pmod 2$. We may assume that $\theta_1\equiv 0\pmod 2$ and $\theta_2\equiv 1\pmod 2$. Then $C_1$ is an even cycle which has exactly two perfect matchings, denoted by $M_{1}$ and $M_{2}$. 
 Let  $M$ be the perfect matching of $C_2\cup (P-V(C_1))$.  Then $\mathbf B(\theta_1,\theta_2;\theta)$ has exactly three perfect 2-matchings: $C_1 \cup M$, $M_{1}\cup M$ and $M_{2}\cup M$. Then, by Theorem~\ref{determinant}, 
 \[\setlength{\jot}{5pt} 
\begin{aligned}
\displaystyle \det(\mathbf B(\theta_1,\theta_2;\theta))&=2 (-1)^{1+|E(C_1)|+|M|} + (-1)^{|M_{1}|+|M|}+(-1)^{|M_{2}|+|M|}\\
&=2 (-1)^{1+\theta_1+(\theta_2+\theta-2)/2}+2(-1)^{(\theta_1+\theta_2+\theta-2)/2}\\
&=2 (-1)^{(\theta_1+\theta_2+\theta-2)/2} ((-1)^{1+\theta_1/2} +1).
\end{aligned}
\]
So $\det(\mathbf B(\theta_1,\theta_1;\theta))\equiv 0 \pmod 2$, and $\det(\mathbf B(\theta_1,\theta_2;\theta))=0$ if and only if $\theta_1\equiv 0 \pmod 4$.

In the following, we consider the case that $\theta$ is even. In this case, $\mathbf B(\theta_1,\theta_2;\theta)$ always has a perfect 2-matching. Let $M$ be the perfect matching of $P-V(C_1\cup C_2)$ that could be empty. 

If $\theta_1\equiv \theta_2 \equiv 1 \pmod 2$, then $\mathbf B(\theta_1,\theta_2; \theta)$ has exactly two perfect 2-matching: its unique perfect matching, and $C_1 \cup C_2\cup M$. Then
\[\det(\mathbf B(\theta_1,\theta_2;\theta))=(-1)^{(\theta_1+\theta_2+\theta-2)/2} + 2^2 (-1)^{2+ \theta_1+\theta_2+ (\theta-2)/2}
\in \{\pm 5, \pm 3\}.\]
So $|\det(\mathbf B(\theta_1,\theta_2;\theta))|\ne 0$ or $1$. 

If $\theta_1\not\equiv \theta_2\pmod 2$, we may assume that $\theta_1\equiv 0 \pmod 2$ and $\theta_2\equiv 1\pmod 2$. Then every perfect 2-matching of $\mathbf B(\theta_1,\theta_2;\theta)$ contains $C_2$. The cycle $C_1$ is an even cycle which has exactly two perfect matchings $M_1$ and $M_2$. Let $M$ be the unique prefect matching of $P-(V(C_1\cup C_2))$. Then $\mathbf B(\theta_1,\theta_2;\theta)$ has exactly three perfect 2-matchings: $C_1\cup C_2 \cup M$, $C_2\cup M_1\cup M$ and $C_2\cup M_2\cup M$. So 
 \[
 \setlength{\jot}{5pt} 
\begin{aligned}
\displaystyle \det(\mathbf B(\theta_1,\theta_2;\theta))&=2^2 (-1)^{2+|E(C_1)|+|E(C_2)|+|M|} + 2 (-1)^{1+|E(C_2)|+|E(C_1)|/2|+|M|}\\
&=4 (-1)^{\theta_1+\theta_2+(\theta-2)/2}+2(-1)^{1+\theta_2+(\theta_1+\theta-2)/2}
\in \{\pm 2, \pm 6\}.
\end{aligned}
\]
It follows that $|\det(\mathbf B(\theta_1,\theta_2;\theta))|\notin \{0, 1\}$.

The last case to consider is that $\theta_1\equiv \theta_2\equiv 0 \pmod 2$. Then both $C_1$ and $C_2$ have exactly two perfect matchings. Let $M$ be the unique perfect matching of $P-V(C_1\cup C_2)$. It follows from $\theta$ is even that $\mathbf B(\theta_1,\theta_2;\theta)$ has exactly nine perfect 2-matchings: one with two cycles together with $M$, two perfect 2-matchings consisting of $C_1$ and a perfect matching of $C_2$ and $M$, two perfect 2-matchings consisting of $C_2$ and a perfect matching of $C_1$ and $M$, and four distinct perfect matchings. Then, by Theorem~\ref{determinant}, 
 \[
 \setlength{\jot}{5pt} 
\begin{aligned}
\displaystyle \det(\mathbf B(\theta_1,\theta_2;\theta))&=2^2 (-1)^{2+|E(C_i)|+|E(C_j)|+|M|} + 2\times 2(-1)^{1+|E(C_i)| +|E(C_j)|/2+|M|}\\
&\, \quad +2\times 2(-1)^{1+|E(C_j)| +|E(C_i)|/2+|M|}+4 (-1)^{(\theta_1+\theta_2+\theta-2)/2}\\
&=4((-1)^{\theta_1+\theta_2+(\theta-2)/2}+(-1)^{\theta_1+(\theta_2+\theta)/2}+(-1)^{\theta_2+(\theta_1+\theta)/2}+ (-1)^{(\theta_1+\theta_2+\theta-2)/2})\\
&= 4(-1)^{(\theta_1+\theta_2+\theta-2)/2} ((-1)^{(\theta_1+\theta_2)/2}+(-1)^{(\theta_1+2)/2}+(-1)^{(\theta_2+2)/2}+1). 
\end{aligned}
\]
So $\det(\mathbf B(\theta_1,\theta_2;\theta))\equiv 0\pmod 4$. Furthermore, $\det(\mathbf B(\theta_1,\theta_2;\theta))=0$ if and only if  $\theta_1$ and $\theta_2$ are not both congruent to 2 modulo 4. This completes the proof.
\end{proof}

\section{Invertibility and sign-invertibility of graphs}

An $(x,y)$-path is a path joining vertices $x$ and $y$. An
$(x,y)$-path $P$ is {\em feasible} if $G-V(P)$ has a perfect 2-matching.  Denote the set of all feasible $(x,y)$-paths of $G$ by $\mathcal{P}_{xy}$. Note that,
if $x=y$, then an $(x,x)$-path is a single vertex $x$ and $\mathcal P_{xx}=\{ x \}$.
 A characterization of the inverse of a graph based on feasible paths and perfect 2-matchings has been given by Yang and the second author in~\cite{YY18} as follows. A general version for weighted graphs can be found in~\cite{YYMK}. 

\begin{thm}[\cite{YY18}]\label{thm:inverse}
Let $G$ be an invertible graph and $A$ be the adjacent matrix of $G$. Then, for any two vertices $x$ and $y$ (they may be the same vertex), 
\[
(A^{-1})_{xy}= 
\displaystyle \frac{1} {\det(G)} \sum_{P\in \mathcal P_{xy}} (-1)^{|E(P)|} \det(G-V(P))  
\]
where $\mathcal P_{xy}$ is the set of all feasible $(x,y)$-paths.
\end{thm}

It follows from the definition of sign-invertibility that $(A^{-1})_{xy}\in \{-1,0,1\}$ for any pair vertices $x$ and $y$ of $G$ if $G$ is sign-invertible. In the rest of this paper, we always use $A$ to denote the adjacency matrix of a graph and $A^{-1}$ to denote the inverse of $A$ without specific mentioning. 

The following result characterizes sign-invertible graphs based on feasible paths.

\begin{thm}\label{thm:s-inverse}
A graph $G$ is sign-invertible if and only if $|\det(G)|=1$, and for any two vertices $x$ and $y$ (they may be the same vertex), 
\[\Big |\sum_{P\in \mathcal P_{xy}} (-1)^{|E(P)|} \det(G-V(P))\Big | \in \big \{ 0,  1 \big \} \]
where $\mathcal P_{xy}$ is the set of all feasible $(x,y)$-paths.
\end{thm}
\begin{proof}
The sufficiency of the theorem follows directly from Theorem~\ref{thm:inverse}. To see the necessity, 
let $G$ be an $s$-invertible graph and $A$ be its adjacency matrix. Then $(A^{-1})_{xy}\in \{-1, 0, 1\}$. Therefore, both $\det(G)=\det(A)$ and $\det(G^{-1})=\det(A^{-1})$ are integer. Note that $\det(A)\det(A^{-1})=1$. So $\det(G)=\det(A)\in\{-1, 1\}$. It follows from Theorem~\ref{thm:inverse} that
\[\Big |\sum_{P\in \mathcal P_{xy}} (-1)^{|E(P)|} \det(G-V(P))\Big |=|(A^{-1})_{xy}| \in \{0, 1\},\]
where $\mathcal P_{xy}$ is the set of all feasible $(x,y)$-paths. This completes the proof.
\end{proof}

Theorem~\ref{thm:s-inverse} shows that every sign-invertible graph is unimodular. It follows from Theorem~\ref{determinant} that a unimodular graph always has a perfect matching. So the following corollary holds. 

\begin{cor}
Every sign-invertible graph is unimodular and hence has a perfect matching. 
\end{cor}

The invertibility and sign-invertiblity of cycles follows directly from  Proposition~\ref{prop:cycle} and Theorem~\ref{thm:s-inverse}.

\begin{prop}\label{prop:c-s-inverse}
A cycle $C_n$ is invertible if and only if $n\not\equiv 0 \pmod 4$. Furthermore, a cycle $C_n$ is not sign-invertible for any $n\ge 3$. 
\end{prop}

The following proposition characterizes the sign-invertibility of theta graphs. 

\begin{prop}\label{prop:theta-s-inverse}
A theta graph $ \Theta(\theta_1, \theta_2, \theta_3)$ is sign-invertible if and only if one of $\theta_i$s is equal to $2$ and the other two are congruent to $0$ modulo $4$. 
\end{prop}
\begin{proof} $\Rightarrow$: Let $G$ be an sign-invertible theta graph. Then $G$ is unimodular.  By Proposition~\ref{prop:theta},  $G=\Theta(\theta_1,\theta_2, \theta_3)$ where $\theta_1, \theta_2$ and $\theta_3$ are even but not all equivalent modulo four. Without loss of generality, assume that $\theta_1\equiv \theta_3\pmod 4$. Let $x$ and $y$ be the two central vertices of $\Theta(\theta_1,\theta_2,\theta_3)$ and $P_1, P_2, P_3$ are the three central paths. 
 \medskip

\noindent{\bf Claim.} {\sl The path $P_2$ has exactly two vertices $x$ and $y$.}\medskip

\noindent{\em Proof of Claim.} Suppose to the contrary that $P_2$ has other two vertices $x'$ and $y'$ such that $x'$ and $y'$ are adjacent to $x$ and $y$ respectively in $P_2$. Note that $\Theta(\theta_1, \theta_2, \theta_3)$ has three feasible $(x',y')$-paths, $Q_1=P_2-\{x,y\}, Q_2=P_1\cup \{xx', yy'\}$ and $Q_3=P_3\cup \{xx',yy'\}$. Then $G-V(Q_1)$ is a cycle of $P_1\cup P_3$ whose length is congruent to $2$ modulo $4$ because $|V(P_1)|\equiv |V(P_3)|\pmod 4$. Then $|\det(G-V(Q_1))|=4$ by Proposition~\ref{prop:cycle}. On the other hand, both $G-V(Q_2)$ and $G-V(Q_3)$ have a perfect matching but no cycles, which implies that $|\det(G-V(Q_2))|=1$ and $|\det(G-V(Q_3))|=1$ by Proposition~\ref{prop:forest}. Therefore, 
\[\Big |\sum_{Q_i\in \mathcal P_{x'y'}} (-1)^{|E(Q_i) |} \det (G-V(Q_i))\Big | \notin \{0,1\}.\]
By Theorem~\ref{thm:s-inverse}, $\Theta(\theta_1, \theta_2, \theta_3)$ is not sign-invertible, which gives a desired contradiction. So the claim follows. 
\medskip

It follows from the claim and Proposition~\ref{prop:theta} that $\theta_1\equiv \theta_3\equiv 0 \pmod 4$. The necessity of the proposition follows. \medskip

\noindent $\Leftarrow$: Assume that $G=\Theta(\theta_1, \theta_2, \theta_3)$ satisfies that one of $\theta_i$s is equal to $2$ and the other two are congruent to $0$ modulo $4$. Without loss of generality, assume that $\theta_2=2$ and $\theta_1\equiv \theta_3\equiv 0 \pmod 4$.  Then $P_2=xy$ and $P_1\cup P_2$ is a cycle of length congruent to 2 modulo 4.  Hence the graph $G$ is bipartite and has an even number of vertices. 

Let $x'$ and $y'$ be two vertices of $G$. Then $x', y' \in P_1\cup P_3$. If $x'$ and $y'$ separates $P_1\cup P_3$ into two paths of odd lengths, then, for any $(x',y')$-path $Q$, the subgraph $G-V(Q)$ has a component with an odd number of vertices and hence has no perfect 2-matchings. Hence any $(x',y')$-path is not feasible. So $\mathcal P_{x'y'}=\emptyset$, and hence 
\[\big |\sum_{Q\in \mathcal P_{x'y'}} (-1)^{|E(Q)|}\det(G-V(Q))\big |=0\] holds automatically. In the following, assume that the two $(x',y')$-paths of $P_1\cup P_3$ have an even number of vertices, which are both feasible. 


First, assume that $G$ has a feasible $(x',y')$-path $Q_1$ such that $G-V(Q_1)$ has a perfect 2-matching with a cycle. Then the cycle is $P_2\cup P_i$ for some $i\in \{1, 3\}$. Without loss of generality, assume $i=3$ (i.e., $Q_1\subseteq P_1$). Note that $|V(P_2\cup P_3)|=|V(P_3)|\equiv 0 \pmod 4$. 
It follows from Proposition~\ref{prop:cycle} that $\det(G-V(Q_1))=0$. Then $G$ has the other two feasible $(x',y')$-paths: $Q_2=P_1\cup P_2- (Q_1-\{x',y'\})$ and $Q_3=P_1\cup P_3-(Q_1-\{x',y'\})$. Then both $G-V(Q_1)$ and $G-V(Q_2)$ has no cycles. Since $|V(P_2)|\not \equiv |V(P_3)|\pmod 4$, it follows that $|V(G-V(Q_2))|/2\not\equiv |V(G)-V(Q_3)|/2 \pmod 2$. Hence $\det(G-V(Q_2))=-\det (G-V(Q_3)) \in \{-1, 1\}$. Therefore, 
\[\big |\sum_{k=1}^3 (-1)^{|E(Q_i)|}\det(G-V(Q_i))\big |=0.\]

In the following, assume that for each feasible $(x',y')$-path $Q$, every perfect 2-matching of $G-V(Q)$ has no cycles. Let $Q_1$ and $Q_3$ be the two feasible $(x',y')$-paths in $P_1\cup P_3$. It follows from $|V(P_1\cup P_3)|\equiv 2\pmod 4$ that $|V(Q_1)| \equiv |V(Q_3)|  \pmod 4$. 

If both $x'$ and $y'$ are contained in the same central path, without loss of generality, assume $x',y'\in P_1$ and $Q_1\subseteq P_1$. Then $G$ has three feasible $(x',y')$-paths: $Q_1$, $Q_3$ and $Q_2=P_1\cup P_2-(Q_1- \{x',y'\})$. Note that $|V(G-V(Q_3))| =|V(Q_1)|-2$ and $|V(G-V(Q_2))|=|V(Q_1)|-2+|V(P_3)|-2\equiv |V(Q_1)|\pmod 4$. It follows that
\[\det(G-V(Q_3))=(-1)^{(|V(Q_1)|-2)/2}= (-1)^{|V(Q_1|/2 +1} =-\det(G-V(Q_2)).\]
Hence 
\[\big |\sum_{k=1}^3 (-1)^{|E(Q_i)|} \det(G-V(Q_i))\big |=|\det(G-V(Q_1))|=1.\] 

Now assume that $x'\in P_1\{x,y\}$ and $y'\in P_3-\{x,y\}$.  Then $x'$ separates $P_1$ into two paths whose lengths have different parities. Without loss of generality, assume that the segment $S_1$ of $P_1$ from $x'$ to $x$ has an odd number of vertices. It follows from that all $(x',y')$-paths on $P_1\cup P_2$ are even that the segment $S_3$ of $P_3$ from $y'$ to $y$ also has an odd number of vertices. So $G$ has exactly three feasible $(x',y')$-paths: $Q_1$, $Q_3$ and $Q_2=S_1\cup P_2\cup S_3$. For convenience, assume $S_1\subset Q_1$ and $S_3\subset Q_3$. Then
$|V(G-V(Q_1))|=(|V(P_1)|-|V(S_1)|-1)+(|V(S_3)|-1)$ and $|V(G-V(Q_2))|=(|V(P_1)|-|V(S_1)|-1)+(|V(P_3)|-|V(S_3)|-1)$. Note that $|V(S_3)|\not\equiv |V(P_3)|-|V(S_3)| \pmod 4$ since $|V(P_3)|=\theta_3\equiv 0\pmod 4$. It follows that
$|V(G-V(Q_1))|\not\equiv |V(G-V(Q_2))|\pmod 4$. Therefore,
\[\det(G-V(Q_1))=-\det(G-V(Q_2)).\]
 Then 
\[\big |\sum_{k=1}^3 (-1)^{|E(Q_i)|} \det(G-V(Q_i))\big |=|\det(G-V(Q_3))|=1.\] 
By Theorem~\ref{thm:s-inverse}, the graph $G$ is sign-invertible and this completes the proof. 
\end{proof}

A connected graph is {\em bicyclic} if it has two edges whose removal results in an acyclic graph. It is not hard to check that a bicyclic graph contains either a barbell or a theta graph as its maximal subgraph without degree-1 vertices. An edge is {\em pendant} if it is incident with a degree-1 vertex. By Propositions~\ref{prop:barbell} and~\ref{prop:theta-s-inverse}, we have the following result on bicyclic graphs.  

\begin{thm}\label{thm:theta-inverse}
A bicyclic graph $G$ without pendant edges is sign-invertible if and only if it is a theta graph $\Theta(2, \theta_2, \theta_3)$ where $\theta_2\equiv \theta_3\pmod 4$. 
\end{thm}

\section{$K_2$-reducible graphs} 

Let $G$ be a graph with a vertex  $v$ of degree one. Assume $u$ is the unique neighbor of $v$ in $G$. A {\em $K_2$-reduction}
of $G$ is the operation of deleting both $u$ and $v$ from $G$. Denote the resulting graph by $G-\{u,v\}$. A subgraph of $G$ is 
{\em $K_2$-irreducible} if it is a minimal subgraph of $G$ obtained from $G$ by a series of $K_2$-reductions. A graph is {\em $K_2$-reducible} if its $K_2$-irreducible subgraph is empty. The following proposition is a direct 
observation.

\begin{prop}\label{prop:reducible}
A $K_2$-reducible graph has a unique perfect $2$-matching that is a perfect matching.
\end{prop}

Let $G$ be a $K_2$-reducible graph. A {\em 2-core} is a maximal subgraph $Q$ of $G$ without vertices of degree one. A connected component of $G-V(Q)$ is called a {\em tree-branch} of $G$, which is a tree. Every tree-branch has exactly one vertex with a unique neighbor in $Q$ that is called {\em the attachment} of the tree-branch. Let $M$ be the unique perfect matching of $G$. For a cycle $C$ of $G$, a {\em peg} of $C$ is an edge of $M$ with exactly one end-vertex on $C$. 


\begin{obs}\label{obs}
Let $G$ be a $K_2$-reducible graph. Then:\\
{\upshape (i)} every cycle of $G$ has at least one peg;\\
{\upshape (ii)} the length of each cycle is congruent to the number of its pegs modulo two;\\
{\upshape (iii)} a tree branch of $G$ has a perfect matching if and only if its attachment is not incident with a peg.
\end{obs}

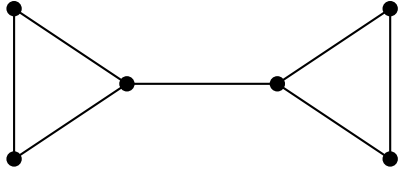
\begin{figure}[htbp]
\begin{center}
\begin{tikzpicture}[thick]
\draw[] (-.5,1)--(-.5,-1)--(1,0)--(-.5,1);
\draw[] (1,0)--(3,0);
 \draw[] (4.5,1)--(4.5,-1)--(3,0)--(4.5,1);

\filldraw[] (-.5,1) circle (2.5pt);
\filldraw[] (4.5,1) circle (2.5pt);
\filldraw[] (-.5,-1) circle (2.5pt);
\filldraw[] (4.5,-1) circle (2.5pt);
\filldraw[] (1,0) circle (2.5pt);
\filldraw[] (3,0) circle (2.5pt);

\end{tikzpicture}
\caption{A $K_2$-irreducible graph with a unique perfect matching.}\label{fig:bowtie}
\end{center}
\end{figure}

There are infinitely many graphs with a unique perfect matching which are not $K_2$-reducible. For example, the graph in Figure~\ref{fig:bowtie}.

\begin{thm}  \label{thm:reduction}
A graph is invertible if and only if it is $K_2$-reducible or its $K_2$-irreducible subgraph is invertible. 
\end{thm}

\begin{proof}
If $G$ is $K_2$-irreducible, then its $K_2$-irreducible subgraph is itself and the theorem holds trivially. So assume that $G$ does has a vertex of degree-1. Let $v$ be a vertex of $G$ with a unique neighbor $u$. Then every perfect 2-matching of $G$ contains the edge $uv$ as a connected component. It follows from Theorem~\ref{determinant} that
\[\det(G)=-\det(G-\{u,v\}).\]
Repeatedly applying the $K_2$-reduction to $G$, let $Q$ be the $K_2$-irreducible subgraph of $G$. If $Q$ is not empty, then it holds that $|\det(G)|=|\det(Q)|$. Hence $G$ is invertible if and only if $Q$ is invertible. If $Q$ is empty, then $G$ has a unique perfect 2-matching which is a perfect matching of $G$. So $|\det(G)|=1$, and hence $G$ is invertible. This completes the proof.
\end{proof}

The above result does not hold for sign-invertibility of graphs. Even a $K_2$-reducible graph may not be sign-invertible. Instead, we have the following weak version of the result.

\begin{thm}\label{thm:unimodular}
Let $G$ be an sign-invertible graph. Then $G$ is $K_2$-reducible or its $K_2$-irreducible subgraph is unimodular.
\end{thm}
\begin{proof}
Let $G$ be an sign-invertible graph. A similar argument as in the above theorem shows that
\[|\det(G)|=|\det(Q)|\]
where $Q$ is the $K_2$-irreducible subgraph of $G$. It follows from Theorem~\ref{thm:s-inverse} that $|\det(G)|=1$. Hence, $Q$ is unimodular. 
\end{proof}


Let $G$ be a $K_2$-reducible graph. A {\em semi-2-core} of $G$ is the subgraph consisting of its 2-core together with all pegs.

\begin{thm}\label{thm:reduce}
Let $G$ be a $K_2$-reducible graph. Then $G$ is sign-invertible if and only if its semi-2-core is sign-invertible. 
\end{thm}
\begin{proof}
Let $G$ be a $K_2$-reducible graph and $M$ be the unique perfect matching of $G$. The result holds trivially if $G$ itself is a semi-2-core. So we may assume that $G$ has a vertex $u$ of degree one such that its neighbor $v$ does not belong to the 2-core of $G$. Then $uv\in M$. Then  it sufficient to prove the following claim.  \medskip

\noindent{\bf Claim:} {\sl the graph $G$ is sign-invertible if and only if that $G-\{u,v\}$ is sign-invertible.} \medskip

\noindent{\em Proof of Claim.} Note that, by Theorem~\ref{determinant}, 
$\det(G)=-\det(G-\{u,v\}$ because every 2-matching of $G-\{u,v\}$ together with the edge $uv$ is a 2-matching of $G$, and vice versa. 

For any two vertices $x$ and $y$ of $G$, if both $x$ and $y$ belongs to $G-\{u,v\}$, then a feasible  $(x,y)$-path of $G$ is also a feasible  $(x,y)$-path of $G-\{u,v\}$ and vice versa. It follows from Theorem~\ref{thm:inverse} that 
\[(A(G)^{-1})_{xy}=(A(G-\{u,v\})^{-1})_{xy}.\] If both $x,y\notin V(G-\{u,v\})$, then $xy=uv$. Then the edge $uv$ is the unique feasible  $(x,y)$-path of $G$. It also follows from Theorem~\ref{thm:inverse} that \[|(A(G)^{-1})_{xy}|=|\det(G-\{x,y\})|.\] By Theorem~\ref{thm:s-inverse}, we only need to consider that exactly one of $x$ and $y$ belongs to $\{u, v\}$. 

Let $P$ be an $(x,y)$-path of $G$. If $uv\notin E(P)$, then $v\in \{x,y\}$ and $G-V(P)$ has an isolated vertex $u$.   So $P$ is not a feasible path. It follows that every feasible  $(x,y)$-path contains the edge $uv$ and $u\in \{x,y\}$.  Without loss of generality, assume $x=u$. Since $uv$ is not a peg, all feasible $(x,y)$-paths contains a neighbor $w$ of $v$ that is different from $u$. Therefore, for every feasible $(x,y)$-path $P$ of $G$, $P'=P-\{u,v\}$ is a feasible $(w,y)$-path of $G-\{u,v\}$. So \[ (-1)^{|E(P)|}=(-1)^{|E(P-\{u,v\})|}=(-1)^{|E(P')|}\] and \[\det(G-V(P))=\det((G-\{u,v\}) -V(P')).\] It follows from Theorem~\ref{thm:inverse} that $|(A(G)^{-1})_{xy}|=|(A(G-\{u,v\})^{-1})_{wy}|$. 

Therefore, the claim follows from Theorem~\ref{thm:s-inverse}, and this completes the proof of the theorem.
\end{proof}

Graphs with a unique perfect 2-matching have been characterized by Yang, Mandal, Klein and the second author in~\cite{YYMK} as follows. 

\begin{thm}[\cite{YYMK}]\label{thm:unique}
A graph has a unique perfect 2-matching  if and only if it is $K_2$-reducible or its $K_2$-irreducible subgraph is a union of disjoint odd cycles. 
\end{thm}

If $G$ is a graph with a unique perfect 2-matching, then $|\det(G)|=2^k$ where $k$ is the total number of cycles in the unique perfect 2-matching. It is obvious that every graph with a unique perfect 2-matching is invertible. As a direct corollary of Theorem~\ref{thm:s-inverse}, the following proposition holds.
 
\begin{prop}
Let $G$ be a graph with a unique perfect 2-matching.  If $G$ is sign-invertible, then it is $K_2$-reducible.
\end{prop}

A $K_2$-reducible bipartite graph is a bipartite graph with a unique perfect matching, and vice versa. The family of $K_2$-reducible bipartite graphs receives particular attention due to many interesting combinatorial properties (see~\cite{CD85, YY18}). It has been proven in~\cite{YY18} that an $(x,y)$-path of a $K_2$-reducible bipartite graph $G$ is feasible if and only if it is an $MM$-alternating path where $M$ is the unique perfect matching of $G$, where a path $P$ is  {\em $MM$-alternating} if $P\cap M$ is a perfect matching of $P$. Let $\tau_o(x,y)$ be the number of all feasible $(x,y)$-paths $P$ with $|E(P)\backslash M| \equiv 1\pmod 2$, and let $\tau_e(x,y)$ be the number of all feasible $(x,y)$-paths $P$ with $|E(P)\backslash M|\equiv 0\pmod 2$. The absolute value $|\tau_o(x,y)-\tau_e(x,y)|$ is called the {\em essential index} of $(x,y)$. 

\begin{thm}
Let $G$ be a $K_2$-reducible bipartite graph. Then $G$ is sign-invertible if and only if the essential index of any pair of two vertices is either 0 or 1. 
\end{thm}
 \begin{proof}
 Let $G$ be a $K_2$-reducible bipartite graph, and let $M$ be its unique perfect matching. So $|\det(G)|=1$.
 
For any two vertices $x$ and $y$, let $P$ be a feasible $(x,y)$-path. Then $P$ is an $M$-alternating and hence $G-V(P)$ is $K_2$-reducible. So
\[\det(G-V(P))= (-1)^{|M|-|E(P)\cap M|}
\]
Therefore, it holds that
\begin{equation*}
\begin{aligned}
\Big |\sum_{P\in \mathcal P_{xy}}(-1)^{|E(P)|} \det(G-V(P))\Big | & = \Big | \sum_{P\in \mathcal P_{xy}} (-1)^{|M|+|E(P)\backslash M|} \Big| \\
& =\Big |\sum_{P\in \mathcal P_{xy}} (-1)^{|E(P)\backslash M|}\Big | \\
&= \Big | \tau_o(x,y)-\tau_e(x,y) \Big |. 
\end{aligned}
\end{equation*}
It follows from Theorem~\ref{thm:s-inverse} that $G$ is sign-invertible if and only if $|\tau_o(x,y) -\tau_e(x,y)|\in \{0,1\}$ for each pair of vertices $(x,y)$ of $G$. 
\end{proof}

It is known that a bipartite graph with a unique perfect matching if and only if it is $K_2$-reducible~\cite{Kotzig}. However, a non-bipartite graph with a unique perfect matching may not be $K_2$-irreducible (for example, see Figure~\ref{fig:bowtie}). It is an interesting problem to ask which non-bipartite graphs with a unique perfect matching is sign-invertible. In a recent paper~\cite{KS}, Kalita and Sarma settled this problem for unicyclic graphs.   \medskip

\section{Sign-invertibility of graphs with small cycle rank}

Let $G$ be a connected graph. The {\em cycle rank} of a graph is the smallest number of edges whose removal results in a subgraph without cycles. In this section, we focus on characterizing  invertibility of graphs with cycle rank at most two.

A graph with cycle rank one is also called a {\em unicyclic graph}, and a graph with cycle rank two is called a {\em bicyclic graph}. It is straightforward to see that the 2-core of a unicyclic graph is a cycle. For a bicyclic graph $G$, let $e_1$ and $e_2$ be the two edges whose removal from $G$ results in a spanning tree $T$. Then $T\cup e_i$ has a unique cycle, denoted by $C_i$ for each $i\in [2]$. If $|V(C_1\cap C_2)|\le 1$, then $C_1$ and $C_2$ are joined by a unique path of $T$ and hence the 2-core of $G$ is a barbell graph. If $|V(C_1 \cap C_2)|\ge 2$, then $C_1\cup C_2$ is a theta graph which is the 2-core of $G$. So we have the following trivial observation. 

\begin{prop}\label{obs}
Let $G$ be a connected graph with cycle rank at most two. Then its 2-core is either a cycle, or a barbell graph, or a theta graph. 
\end{prop}

The following proposition provides information of feasible paths in $K_2$-reducible graphs with cycle rank at most two. 

 \begin{prop}\label{prop:3-peg}
Let $G$ be a $K_2$-reducible graph with cycle rank at most two. If every cycle of $G$ has at least three pegs, then, for any vertices $x$ and $y$, there is at most one feasible $(x,y)$-path $P$. Furthermore, for a feasible path $P$, the subgraph $G-V(P)$ is $K_2$-reducible.
\end{prop}
\begin{proof}
Let $G$ be a $K_2$-reducible graph with cycle rank at most two, and let $Q$ be the 2-core of $G$. Suppose to the contrary that there exists two vertices $x$ and $y$ such that $G$ has at least two feasible $(x,y)$-paths $P$ and $P'$.

The symmetric difference $P\Delta P'$ has a cycle $C$. Since $C$ has at least three pegs, it follows that at least one peg  $uv$ of $C$ does not belong to $P$ nor $P'$. Then one end-vertex of $uv$ belongs to either $P$ or $P'$. Without loss of generality, assume that $u\in V(P)$ and $v\notin C$. If $u$ is an attachment of a tree-branch $T$, then $T$ has an odd number of vertices and has no perfect 2-matching, which contradicts that $P$ is a feasible path. So $u$ is not an attachment of a tree-branch. It follows that $uv$ belongs to $Q$. So $Q$ is not a cycle. By Proposition~\ref{obs}, the 2-core $Q$ is either a barbell or a theta graph. 

If $Q$ is a barbell, then $Q$ has another cycle $C'$ joining to $C$ by a path  which contains $uv$. Since $C'$ also has at least three pegs, it follows that the component of $G-V(P)$ containing $C'$ does not have a perfect 2-matching, which contradicts that $P$ is feasible.

So assume that $Q$ is a theta graph, and let $W$ be the path of $Q$ joining two vertices of $C$, one of which is $u$. Note that $uv\in W$. It further follows that all pegs of $C$ not on $P$ or $P'$ belong to $W$. Hence $C$ has at most four pegs (at most two belong to $P\cap P'$ and at most two belong to $W$).
Let $W_1$ and $W_2$ be the two segments of $C$ separated by the two end-vertices of $W$. Since every cycle of $G$ has at least three pegs and all pegs of $C$ belong to $P\cup P' \cup W$, it follows that one of the two cycles $W_1\cup W$ and $W_2\cup W$ has a peg $ww'$ with an end-vertex $w$ inside of $W$ (exclude two end-vertices of $W$). Then the segment $W(u,w)$ of $W$ from $u$ to $w$ (not including the end-vertices) has an odd number of vertices. Note that $G-V(P)$ is acyclic and has no perfect matching due to $W(u,w)$ has no perfect matching, which contradicts that $P$ is a feasible path. The contradiction implies that $G$ has at most one feasible $(x,y)$-path.

Now, assume that $P$ is a feasible path of $G$. For any cycle $C$ of $G$, the subgraph $C-V(P)$ is either acyclic or $P$ does not contain pegs of $C$.  No matter which case happens, $G-V(P)$ is $K_2$-reducible because $P$ is feasible. This completes the proof.
\end{proof}

The following result provides a partial result on sign-invertible graphs with small cycle rank. 

\begin{thm}\label{thm:3-peg}
Let $G$ be a $K_2$-reducible graph with cycle rank at most two. If every cycle of $G$ has at least three pegs, then $G$ is sign-invertible. 
\end{thm}
\begin{proof}
Let $G$ be a $K_2$-reducible graph with cycle rank at most two, and let $x$ and $y$ be two vertices of $G$. Then it follows from Proposition~\ref{prop:3-peg} that $G$ has at most one feasible $(x,y)$-path $P$, and furthermore, if such feasible path $P$ exists, then $G-V(P)$ is $K_2$-reducible. Hence, either $(A^{-1})_{xy}=0$ or
\[(A^{-1})_{xy}=(-1)^{|EP|} (-1)^{|V(G)-V(P)|/2} \in \{-1,1\}.\]
By Theorem~\ref{thm:s-inverse}, the graph $G$ is sign-invertible. 
\end{proof}

In the following, we focus on  invertibility of unicyclic graphs and bicyclic graphs.

\subsection{Unicyclic graphs}

The invertibility of unicyclic graphs is given in the following result. 
 
 \begin{thm}\label{thm:uni-inverse}
A unicyclic graph $G$ is invertible if and only if $G$ is $K_2$-reducible or its $K_2$-irreducible subgraph is not a $4k$-cycle. 
 \end{thm}
\begin{proof} Let $G$ be a unicyclic graph. 

First, assume that $G$ is invertible. If $G$ is not a $K_2$-reducible graph, then let $H$ be its $K_2$-irreducible subgraph. Then $H$ is a cycle, the unique cycle of $G$. It follows from Proposition~\ref{prop:cycle} that $H$ is not a $4k$-cycle.

Now assume that $G$ is $K_2$-reducible or the $K_2$-irreducible subgraph of $G$ is not a $4k$-cycle. If $G$ is $K_2$-reducible, Theorem~\ref{thm:reduction} implies that $G$ is invertible. If $G$ is not $K_2$-reducible, let $H$ be its $K_2$-irreducible subgraph. Then $H$ is a cycle with length not divisible by four. It follows from Theorem~\ref{thm:reduction} and Proposition~\ref{prop:cycle} that $G$ is invertible. This completes the proof.          
\end{proof}


In the rest of this subsection, we focus on sign-invertibility of unicyclic graphs. If $G$ is $K_2$-reducible unicyclic graph, then it has a unique perfect matching $M$ by Proposition~\ref{prop:reducible}. Recall that a peg is an edge of $M$ which  joins a vertex of a tree branch and its attachment, and a path $P$ is $MM$-alternating if $P\cap M$ is a perfect matching of $P$.

\begin{prop}\label{prop:feasible}
Let $G$ be a $K_2$-reducible unicyclic graph with at least two pegs, and let $M$ be the unique perfect matching of $G$. Then,\\
{\upshape(i)} for any vertex $x$, the subgraph $G-x$ has no perfect 2-matchings;\\ 
{\upshape (ii)} for any two different vertices $x$ and $y$, a feasible $(x,y)$-path $P$ is $MM$-alternating and $G-V(P)$ is $K_2$-reducible. 
\end{prop}
\begin{proof}
Let $G$ be a $K_2$-reducible unicyclic graph with at least two pegs $u_1v_1$ and $u_2v_2$ such that $v_1$ and $v_2$ are on the unique cycle $Q$ of $G$. Let $T_k$ be the tree branch of $G$ containing $u_k$ for $k\in [2]$. Let $x$ be a vertex of $G$. Then one of $T_1\cup u_1v_1$ and $T_2\cup u_2v_2$ does not contain the vertex $x$, say $T_1\cup u_1v_1$. Let $M_1=M\cap (T_1\cup u_1v_1)$. Then $G-V(M_1)$ is acyclic. It follows that $G-V(M_1)-x$ has no perfect 2-matchings. Hence $G-x$ does not have a perfect 2-matching and (i) follows.

For any two distinct vertices $x$ and $y$, let $P$ be a feasible $(x,y)$-path. Then, for each peg of $G$, the path $P$ either contain the peg or is disjoint from the peg. Since each peg is a cut edge, $P$ contains at most two pegs. 

If $P$ does not contain all pegs, then the existence of a peg not in $P$ implies that $G-V(P)$ is a forest with a perfect 2-matching because $P$ is feasible. So $G-V(P)$ is $K_2$-reducible and hence has a unique perfect matching $M'= M \cap G-V(P)$. Then $M-M'=M\cap P$, and hence $P$ is $MM$-alternating. Therefore the proposition holds. 

So assume that $G$ has exactly two pegs and $P$ contains both of $u_1v_1$ and $u_2v_2$. Then $T_k-V(P)$ has a perfect matching for both $k\in [2]$ since $P$ is a feasible $(x,y)$-path. Since $P\cap Q\ne \emptyset$, it follows that $G-V(P)$ is a forest with perfect 2-matching. So $G-V(P)$ is $K_2$-reducible and has a perfect matching $M'=M\cap (G-V(P))$. Further, $M-M'=M\cap P$ is a perfect matching of $P$, which implies that $P$ is $MM$-alternating. This completes the proof.
\end{proof}

The above proposition does not hold for $K_2$-reducible unicyclic graph with exactly one peg. For example, the unicyclic graph obtained by joining an isolated vertex to a vertex of an odd cycle.  

\begin{thm}\label{thm:unicycle}
Let $G$ be a unicyclic graph and let $Q$ be the unique cycle of $G$. Then $G$ is sign-invertible if and only if the number of pegs of $Q$ is at least two pegs and the equality holds only if the attachments of two pegs separates $Q$ into two paths of even lengths not congruent modulo $4$. 
\end{thm}
\begin{proof}
Let $G$ be an sign-invertible unicyclic graph. By Proposition~\ref{prop:c-s-inverse}, $G$ is not a cycle. If $G$ is not $K_2$-reducible, then its $K_2$-irreducible subgraph is $Q$ because $G$ is unicyclic. It follows from Theorem~\ref{thm:unimodular} and Proposition~\ref{prop:cycle} that $G$ is not unimodular and hence not sign-invertible.

In the following, we assume that $G$ is $K_2$-reducible. 
By Proposition~\ref{prop:reducible}, $G$ has a unique perfect 2-matching which is a perfect matching $M$ of $G$. Then $|\det(G)|=|(-1)^{|M|}|=1$ by Theorem~\ref{determinant}. By Theorem~\ref{thm:reduce}, assume that $G$ consists of $Q$ together with some pendent edges which are incident with exactly one vertex of $Q$, which are pegs of $Q$.

First, assume that $Q$ has exactly one peg $uv$ such that $v$ is on the unique cycle $Q$ of $G$. Then $Q$ is an odd cycle. Further, $G-u$ has a perfect 2-matching consisting of a cycle $Q$ and a matching.
Hence
\[\Big |\sum_{P\in \mathcal P_{uu}}(-1)^{|E(P)|}\det(G-V(P))\Big |=\Big |\det(G-u)\Big |=\big | \det(Q) \big |=2 \notin \{0,1\}\] 
contradicting Theorem~\ref{thm:s-inverse}. 

By Theorem~\ref{thm:3-peg}, assume that  $G$ has exactly two pegs $u_1v_1$ and $u_2v_2$ with $v_1$ and $v_2$ on $Q$. 
By Proposition~\ref{prop:feasible}, for any feasible $(x,y)$-path $P$, the subgraph $G-V(P)$ is $K_2$-reducible and such path does not exist if $x=y$. So $|\sum_{P\in\mathcal P_{xx}} \det(G-V(P))|=0$. So assume that $x\ne y$. Since $G$ is a unicyclic graph, $G$ has at most two feasible $(x,y)$-paths. Then, for each feasible $(x,y)$-path $P$, it always holds that $|\det(G-V(P))|=1$ because $G-V(P)$ is $K_2$-reducible. If $G$ has exactly two feasible $(x,y)$-paths $P_1$ and $P_2$, then the symmetric difference $P_1\Delta P_2=Q$ and both $u_1v_1$ and $u_2v_2$ belong to $P_1$ and $P_2$ by Proposition~\ref{prop:feasible}. The feasibility of $P_1$ and $P_2$ implies that $|V(P_1\cap Q)|\equiv |V(P_2\cap Q)|\equiv 0 \pmod 2$, and hence $|E(P_1)|\equiv |E(P_2)|\pmod 2$. It follows from Theorem~\ref{thm:s-inverse} that $G$ is sign-invertible if and only if 
\[(-1)^{|E(P_1)|}\det(G-V(P_1))\ne (-1)^{|E(P_2)|} \det(G-V(P_2)),\] equivalently, 
\[(-1)^{|M|- |V(P_1)|/2} \ne (-1)^{|M|-|V(P_2)|/2}.\] 
Hence, $G$ is sign-invertible if and only if $|V(P_1)|\not\equiv |V(P_2)| \pmod 4$ (or equivalently $|V(P_1\cap Q)|\not\equiv |V(P_2\cap Q)|\pmod 4$). In other words, if $G$ has exactly two pegs, then $G$ is sign-invertible if and only if the numbers of vertices of the two $(v_1,v_2)$-paths of $Q$ are both even but not congruent modulo 4. This completes the proof. 
\end{proof}

The following result is a direct corollary of Theorem~\ref{thm:unicycle} and Observation~\ref{obs}.

\begin{cor}[Kalita and Sarma, \cite{KS}]
Let $G$ be a non-bipartite unicylic graph with a unique perfect matching. Then $G$ is sign-invertible if and only if it has at least three pegs. 
\end{cor}

\subsection{Bicyclic graphs}

A bicyclic graph $G$ is a connected graph with cycle rank two, i.e., a bicyclic graph could be  reduced to a tree by deleting exactly two edges. 
 
\medskip

\noindent{\bf Definition.} A cycle $C$ satisfies the {\bf peg condition} if $C$ has either three pegs or two pegs which separates $C$ into two paths of even lengths not congruent to each other modulo four. \medskip

So a unicyclic graph is sign-invertible if and only if it satisfies the peg condition. The similar result could be proved for bicyclic graphs.


\begin{thm}\label{lem:s-inv-barbell}
Let $G$ be a bicyclic graph with barbell as its 2-core. Then $G$ is sign-invertible if and only if $G$ is $K_2$-reducible and satisfies one of the following:\\
{\upshape (i)}  the graph $G$ is obtained by attaching one pedant edge to the central vertex of  $\mathbf B(\theta_1,\theta_2;1)$ where  $\theta_1$ and $\theta_2  \pmod 4$ are odd but not congruent modular four; \\
{\upshape (ii)} every cycle satisfies the peg condition.
\end{thm}
\begin{proof} $\Rightarrow$: 
Assume that $G$ is sign-invertible and let $Q$ be its 2-core, a barbell $B(\theta_1, \theta_2; \theta)$. By Propositions~
 \ref{prop:cycle}, \ref{prop:barbell}, and Theorem~\ref{thm:unimodular}, the graph $G$ is $K_2$-reducible. So every cycle of $G$ has at least one peg. By Theorem~\ref{thm:reduce}, assume $G$ itself is a semi-2-core. Then every tree-branch of $G$ is a single vertex. 
 
If $G$ has a cycle $C$ with exactly one peg $uv$ such that $v\in C$, then both $C$ and $G-u$ have an odd number of vertices. Then $G-u$ has no perfect matching and is not $K_2$-reducible. Then the irreducible subgraph of $G-u$ is either a barbell or the cycle $C$. If the irreducible subgraph of $G-u$ is the cycle $C$, which has an odd number of vertices, then it follows from Theorem~\ref{thm:s-inverse} and Proposition~\ref{prop:cycle} that $(A^{-1})_{uu} \in \{-2, 2\}$, a contradiction to that $G$ is sign-invertible. So assume that the irreducible graph of $G-u$ is the barbell $B(\theta_1, \theta_2; \theta)$ with an odd number of vertices, which implies that $\det(B(\theta_1, \theta_2; \theta))=0$ by Theorem~\ref{determinant} and Theorem~\ref{thm:s-inverse}. Since $C$ is an odd cycle, without loss of generality, assume $\theta_1=|V(C)|\equiv 1\pmod 2$. Then it follows from (ii) and (iii) of Proposition~\ref{prop:barbell} that $\theta\equiv 1\pmod 2$. If $\theta\ge 3$, let $P$ be the central path. Since $G$ is $K_2$-reducible and $C$ has only one peg $uv$, it follows that $P$ contains a peg of another cycle $C'$ of $B(\theta_1, \theta_2; \theta)$. Note that $B(\theta_1, \theta_2; \theta)$ is the $K_2$-irreducible graph of $G-u$. So $C'$ has exactly one peg which belongs to $P$. Therefore, $G= B(\theta_1, \theta_2; \theta) \cup uv$. Let $u'v'$ be the peg of $C'$ such that $v'\in C'$. Then $G-u'$ has a unique perfect 2-matching which has an odd cycle $C'$ as a component. It follows that $(A^{-1})_{u'u'}\in \{-2, 2\}$, a contradiction to that $G$ is sign-invertible. So it follows that $\theta=1$. Then both cycles $C$ and $C'$ of $B(\theta_1, \theta_2; \theta)$ are odd cycles. By (ii) of Proposition~\ref{prop:barbell}, it holds that $\theta_1\not\equiv \theta_2 \pmod 4$. Let $v'$ be the central vertex of $B(\theta_1, \theta_2; 1)$ and $u'v'\in M$. If $v\ne v'$, let $P$ be an $MM$-alternating $(u,v')$-path. Then $P\subseteq C\cup uv$, and $P'=P-w$ is a feasible $(u,u')$-path of $G$. Note that $G-V(P')$ has a unique perfect 2-matching with $C'$ as a component. It follows $|(A^{-1})_{uu'}|=2$ by Proposition~\ref{prop:cycle}, a contradiction to the sign-invertibility of $G$ again. The contradiction implies that $v=v'$ and (i) follows.

Now assume that every cycle of $G$ has at least two pegs.  If $G$ has a cycle $C$ with exactly two pegs $uv$ and $u'v'$ such that the two attachments, say $v$ and $v'$, separate $C$ into two paths  of even lengths congruent modulo four, then $G$ has exactly two feasible $(u,u')$-paths $P_1$ and $P_2$ such that $E(P_1)\Delta E(P_2)=E(C)$ and hence $|E(P_1)|\equiv |E(P_2)| \pmod 4$. Therefore, 
\begin{equation*}
\begin{aligned}
\det(G-V(P_1)) (-1)^{|E(P_1)|} &=\det\big(G-(V(C)\cup \{u, u'\})\big) (-1)^{(|V(P_2)|-4)/2} (-1)^{|E(P_1)|}\\
&=\det\big(G-(V(C)\cup \{u, u'\})\big) (-1)^{(|E(P_2)|-3)/2+|E(P_1)|}\\
&= \det\big(G-(V(C)\cup \{u, u'\})\big) (-1)^{(|E(P_1)|-3)/2+|E(P_2)|}\\
&=\det(G-V(P_2))(-1)^{|E(P_2)|} \ne 0.
\end{aligned}
\end{equation*}
It follows from Theorem~\ref{thm:inverse} that
\[|(A^{-1})_{uu'}|=\sum_{i=1}^2\det(G-V(P_i)) (-1)^{|E(P_i)|}=2\det(G-V(P_1)) (-1)^{|E(P_1)|} \notin\{-1,0,1\},\]
which contradicts the sign-invertibility of $G$ and Theorem~\ref{thm:s-inverse}. So every cycle of $G$ satisfies the peg condition, and (ii) holds.  \medskip

\noindent $\Leftarrow$: Let $G$ be a bicyclic graph with a barbell as its 2-core $Q$. Assume that $G$ is $K_2$-reducible and satisfies either (i) or (ii). Let $M$ be the unique perfect matching of $G$. By Theorem~\ref{thm:reduce}, we may assume that $G$ consists of its 2-core $Q$ together with all pegs of cycles of $G$. By Theorem~\ref{thm:s-inverse}, it suffices to verify that $(A^{-1})_{xy}\in \{-1,0,1\}$ for any two vertices $x$ and $y$ of $G$. Let $P$ be an $(x,y)$-path (it may holds that $x=y$). 

First, assume (i) holds. Then  $Q= \mathbf B(\theta_1, \theta_2;1)$. Let $v$ be the central vertex and $u$ be the vertex of degree one. Then $uv\in M$. If $P$ does contain the central vertex $v$, then $G-V(P)$ is acyclic. Hence $P$ is feasible if and only if it is an $MM$-alternating. If $P$ does not contain $u$ and $v$, then $P$ is feasible if and only if $G-V(P)$ is $K_2$-reducible. Hence, if $P$ is not a single vertex $u$, then $P$ is feasible if and only if it is $MM$-alternating. Hence $G$ has at most one feasible $(x,y)$-path $P$ and $G-V(P)$ has a unique perfect 2-matching which is a matching. It follows that $(A^{-1})_{xy}\in \{-1, 0,1\}$. If $P=u$ (i.e. $x=y=u$), then $G-x=\mathbf B(\theta_1,\theta_2;1)$. By (ii) of Proposition~\ref{prop:barbell}, it holds that $(A^{-1})_{xx}=0$. So $G$ is sign-invertible by Theorem~\ref{thm:s-inverse}. 

In the following, assume (ii) holds. The peg condition implies that every cycle of $G$ has at least two pegs. By Proposition~\ref{prop:3-peg} and Theorem~\ref{thm:s-inverse}, it suffices to prove the theorem holds for bicyclic graphs  having one cycle $C$ with exactly two pegs. 

If $G$ has no feasible $(x,y)$-path, it holds trivially that $(A^{-1})_{xy}=0$. Assume that $P$ is a feasible $(x,y)$-path.  For a peg of some cycle of $G$, the path $P$ either contains the peg or is disjoint from the peg. If every cycle of $G$ has a peg disjoint from $P$, then $P$ is the unique feasible $(x,y)$-path and $G-V(P)$ is $K_2$-reducible. So $(A^{-1})_{xy} \in \{-1,1\}$. Therefore, assume that $G$ has a cycle such that $P$ does contain all pegs of the cycle. Since $G$ has barbell as its 2-core, it follows that the cycle with all pegs in $P$ has exactly two pegs. Without loss of generality, assume the cycle is $C$. Then $P'=P\Delta E(C)$ is another feasible $(x,y)$-path. Furthermore, both $P$ and $P'$ are $MM$-alternating paths. By the peg condition, $|E(P')|\not\equiv |E(P)|\pmod 4$. So, by Theorem~\ref{determinant}, 
\[(-1)^{|E(P)|}\det(G-V(P)) =
(-1)^{|M|+|E(P)|/2}=(-1) ((-1)^{|E(P')|}\det(G-V(P'))).\]
It follows from Theorem~\ref{thm:s-inverse} that $(A^{-1})_{xy}=0$.   This completes the proof.
\end{proof} 

In the following, we consider the last case -- bicyclic graphs with theta graphs as its 2-core. \medskip

\noindent{\bf Definition.} A $K_2$-reducible bicyclic graph $G$ with a theta graph $\Theta(\theta_1, \theta_2, \theta_3)$ as 2-core is called a {\bf crab-graph} if:\\ 
(i) $G$ has exactly two vertices of degree one $x$ and $y$; and \\
(ii) $G$ has exactly three feasible $(x,y)$-paths which all have an even number of vertices but not congruent to each other modular four. \medskip

It is not hard to see that a crab-graph $G$ has a unique perfect matching $M$ and all three feasible $(x,y)$-paths are $MM$-alternating. Furthermore, two degree-one vertices belong to the same central path of the theta graph of $G$.

\begin{thm}\label{lem:s-inv-theta}
Let $G$ be a bicyclic graph with theta graph $\Theta(\theta_1, \theta_2, \theta_3)$ as its 2-core. Then $G$ is sign-invertible if and only if 
either the $K_2$-irreducible subgraph of $G$ is $\Theta(2,\theta_2,\theta_3)$ with $\theta_2\equiv \theta_3\equiv 0\pmod 4$,  or $G$ is $K_2$-reducible and its semin-2-core satisfies one of the following:\\
{\upshape (i)} the semi-2-core of $G$ has an exactly one degree-1 vertex with a neighbor on an odd central path of $\Theta(\theta_1, \theta_2, \theta_3)$ where $\theta_1\equiv 1 \pmod 2$ and even $\theta_2\not \equiv \theta_3 \pmod 4$; or\\
{\upshape (ii)}  every cycle of $G$ satisfies the peg condition unless the semi-2-core of $G$ is a crab-graph.
\end{thm}

\begin{proof}
By Theorem~\ref{thm:reduce}, we may assume that $G$ consists of a theta graph $Q$ and all pegs. If $G=Q$, then it follows from Theorem~\ref{thm:theta-inverse} that $G$ is sign-invertible if and only if $Q=\Theta(2, \theta_2, \theta_3)$ with $\theta_2\equiv \theta_3 \equiv 0 \pmod 4$. So in the following, assume that $G$ is $K_2$-reducible and then $G$ has a unique perfect matching $M$. Let $P_i$ be the central path of $Q$ with $\theta_i$ vertices for $i\in [3]$, and let $x$ and $y$ be the two central vertices. \medskip

\noindent $\Rightarrow$: Assume that $G$ is sign-invertible. If $G$ has an exactly one vertex of degree-1, say $u$. Let $v$ be the neighbor of $u$. Then $uv\in M$ and $G-u$ has an odd number of vertices and $\det(G-u)\in \{-1, 0, 1\}$. Then one of $\theta_i$s of $\Theta(\theta_1, \theta_2, \theta_3)$ is odd, say $\theta_1$. By Proposition~\ref{prop:theta}, either $\theta_1\equiv \theta_2\equiv \theta_3\equiv 1\pmod 2$, or  $\theta_2$ and $\theta_2$ are even but not congruent modular four. Note that $G$ is $K_2$-reducible and the central path $P_1$ has an odd number of vertices, the cycle $P_2\cup P_3$ has two pegs on $P_1$. It follows that both $P_2$ and $P_3$ have an even number of vertices. So $\theta_2$ and $\theta_3$ are even and $\theta_2\not\equiv \theta_3\pmod 4$. If $v\notin V(P_1)$, then $v$ is an interior vertex of either $P_2$ or $P_3$, say $P_2$. Since $G$ is $K_2$-reducible, the cycle $P_1\cup P_3$ has exactly one peg on $P_2$. Let the peg be $xx'\in M\cap E(P_2)$. Then $G-x'$ has a unique perfect 2-matching with the odd cycle $P_1\cup P_3$ as a component. It follows from Theorem~\ref{thm:inverse} that $(A^{-1})_{x'x'} \in \{-2,2\}$, contradicting the sign-invertibility of $G$ and Theorem~\ref{thm:s-inverse}. The contradiction implies that $v\in V(P_1)$ and (i) holds. 

Now, assume that $G$ has at least two vertices of degree one. If $G$ has a cycle $C$ with exactly one peg $uv$ with $v\in V(C)$, then $C$ is an odd cycle. Then $G-u$ has no perfect matching but has a degree-1 vertex. So $G-v$ does have a unique perfect 2-matching $(M-(M\cap E(C))\cup \{uv\})\cup \{C\}$. Then, by Proposition~\ref{prop:cycle}, 
\[(A^{-1})_{uu}=(-1)^{|E(P)|}\det(G-V(P))\in \{-2, 2\},\] 
contradicting the sign-invertibility of $G$ and Theorem~\ref{thm:s-inverse}. Hence, every cycle of $G$ has at least two pegs. 

By the peg condition, we may assume that $G$ has a cycle $C$ with exactly two pegs $u_1v_1$ and $u_2v_2$ such that $v_1, v_2\in V(C)$. Then $u_1$ and $u_2$ are joined by at least two feasible paths whose symmetric differences is $C$. Note that both the two feasible $(u_1,u_2)$-paths are $MM$-alternating.   

If $G$ has exactly two feasible $(u_1,u_2)$-paths, it follows from the $K_2$-reducibility of $G$ that removing all vertices of each feasible paths results in an acyclic graph. (Otherwise, the perfect 2-matching in the subgraph contains an even cycle, contradicting that $G$ has a unique perfect matching.)
It follows from Theorem~\ref{thm:s-inverse} that $u_1$ and $u_2$ separate $C$ into two paths of even lengths not congruent modulo four. The peg condition holds.

Since $G$ has theta graph as its 2-core, it has at most three feasible $(u_1,u_2)$-paths. So assume that $G$ has exactly three feasible $(u_1,u_2)$-paths. It is not hard to derive that both $v_1$ and $v_2$ belong to the same central path of $\Theta(\theta_1, \theta_2, \theta_3)$. The union of the three feasible paths is the graph $G$, and all interior vertices of the feasible paths belong to the 2-core $Q$. It follows that $u_1$ and $u_2$ are the only degree-1
 vertices. So $G$ consists of $Q=\Theta(\theta_1,\theta_2,\theta_3)$ together with the two pendant edges $v_1u_1$ and $v_2u_2$. 
Note that the three feasible $(u_1,u_2)$-paths are all $MM$-alternating. Let   $\alpha_1, \alpha_2$ and $\alpha_3$ be the number of vertices of the three feasible paths, respectively.  Let $P$ be the feasible $(u_1, u_2)$-path with $\alpha_i$ vertices. Then $G-V(P)$ is $K_2$-reducible because every cycle has at least two pegs. 
\[(-1)^{|E(P)|}\det(G-V(P)) = (-1)^{\alpha_i-1} (-1)^{|M|-\alpha_i/2} =(-1)^{|M|-1+\alpha_i/2}.\] 
So, it follows that
\begin{equation*}
\begin{aligned}
(A^{-1})_{u_1u_2}&=(-1)^{|M|-1+\alpha_1/2} + (-1)^{|M|-1+\alpha_2/2} +(-1)^{|M|-1+\alpha_3/2}\\
&=(-1)^{|M|-1} ((-1)^{\alpha_1/2}+(-1)^{\alpha_2/2}+(-1)^{\alpha_3/2}).
\end{aligned}
\end{equation*}
Hence $(A^{-1})_{u_1u_2}\in \{-1, 0, 1\}$ if and only if $\alpha_1, \alpha_2$ and $\alpha_3$ are even but not all congruent to each other modulo four. Hence (ii) follows.  \medskip


\noindent $\Leftarrow$: By Theorem~\ref{thm:theta-inverse}, it suffices to prove that a $K_2$-reducible graph $G$ satisfying either (i) or (ii) is sign-invertible. In other words, for any two vertices $x$ and $y$ of $G$ (they may be the same vertex), it is sufficient to show $(A^{-1})_{xy}\in \{-1,0,1\}$ by Theorem~\ref{thm:s-inverse}.  By Theorem~\ref{thm:3-peg}, assume that $G$ has a cycle with at most two pegs.

If $G$ has no feasible $(x,y)$-path, then $(A^{-1})_{xy}=0$. So we assume that $G$ has a feasible $(x,y)$-path $P$. For a peg of a cycle of $G$, then either $P$ contains the peg or $P$ is disjoint from the peg. 

First, assume that $P$ is an $MM$-alternating path. Then $G-V(P)$ has a unique perfect matching and at most one cycle because $G$ is $K_2$-reducible and consists of a theta graph and pegs. Further, $G-V(P)$ is also $K_2$-reducible. If $P$ is the unique feasible $(x,y)$-path, then $(A^{-1})_{xy}\in \{-1,1\}$. If $P$ is not the unique feasible $(x,y)$-path, let $P'$ be the other feasible $(x,y)$-path. Then $P\Delta P'$ is a cycle $C$ of $G$ since the 2-core of $G$ is a theta graph. 

If $P'$ is also $MM$-alternating, then both $x$ and $y$ are degree-1 vertices. If $P$ and $P'$ are the only feasible $(x,y)$-paths, it follows from the peg condition in (ii) that 
\[(-1)^{|E(P)|}\det(G-V(P))=(-1)^{|M|}(-1)^{|E(P)|/2}=(-1)^{|M|}(-1)^{|E(P')|/2} (-1)=-(-1)^{|E(P')|}\det(G-V(P')).\]
Hence, it follows from Theorem~\ref{thm:inverse} that $(A^{-1})_{xy}=0$. Since the 2-core of $G$ is a theta graph,  $G$ has at most three feasible $(x,y)$-paths. If $G$ has exactly three feasible $(x,y)$-paths including $P$ and $P'$, then all of them are $MM$-alternating. It follows from (ii) and Theorem~\ref{thm:inverse} that $(A^{-1})_{xy}\in \{-1,1\}$. 

So in the following assume that $P'$ is not $MM$-alternating. Then $G-V(P')$ has no perfect matching. Hence a perfect 2-matching of $G-V(P')$ contains exactly one odd cycle $C'$, the only cycle disjoint from $P'$, which implies that $G-V(P')$ has a unique perfect 2-matching. So $C'$ does not contain attachments of any tree-branches. Note that, all interior vertices of $P$ and $P'$ are not attachments of tree-branches. Hence $G$ has at most two tree-branches, which are $x$ and $y$. If both $x$ and $y$ are degree-1 vertices, then both $P$ and $P'$ are $MM$-alternating, a contradiction to that $P'$ is not $MM$-alternating. Hence exactly one of $x$ and $y$ has degree one, which contradicts (i) and that $C'$ is an odd cycle of $G-V(P')$. 

So assume that all feasible $(x,y)$-paths $P$ are not $MM$-alternating, which including the case $x=y$. Since $G$ has a unique perfect matching $M$ and $P$ is not $MM$-alternating, it follows that every perfect 2-matching of $G-V(P)$ contains exactly one  odd cycle and $P$ has an odd number of vertices. Since the 2-core of $G$ is a theta graph, $G$ has at most two different odd cycles and hence $G-V(P)$ has at most two perfect 2-matchings. If $G-V(P)$ does have exactly two perfect 2-matchings, then the union of the two odd cycles of the perfect 2-matching is the 2-core of $G$. So $P$ has to be a single vertex $x$ (i.e, $x=y$). Then $G$ has exactly one vertex of degree one. By (i) and Proposition~\ref{prop:barbell}, it holds that $(A^{-1})_{xx}\in \{-1,1\}$. If $G-V(P)$ has exactly one perfect 2-matching with an odd cycle, which has exactly one peg incident with one of the central vertex of the 2-core $Q$ of $G$, then $G$ does not satisfies (i) and the odd cycle in the perfect 2-matching has exactly one peg. So the odd cycle also does not satisfy the peg condition. Hence $G$ does not satisfies (ii), a contradiction. This completes the proof. 
\end{proof}

A characterization of sign-invertible bicyclic graphs follows directly from combining Theorem~\ref{lem:s-inv-barbell} and Theorem~\ref{lem:s-inv-theta}. \medskip

\noindent{\bf Remark.} The characterization of graphs with cycle rank at most two involves quite heavy structure analysis. The $K_2$-reduction makes the structures of these graphs much simpler. The sign-invertibility for graphs with high cycle ranks becomes more complicated. It is interesting to study the sign-invertibility of graphs with high cycle ranks. The family of cactus may be particularly interesting due to their structures, where a {\em cactus} is a connected graph in which any two simple graphs have at most one vertex in common. The methods developed in this paper may be useful to study the sign-invertibility of these graphs.

\end{document}